\documentclass[twoside,12pt]{article}

\usepackage{amsmath}
\usepackage{amssymb}
\usepackage{amscd}
\usepackage{amsthm}

\numberwithin{equation}{subsection}

\theoremstyle{plain}

\theoremstyle{remark}

\theoremstyle{definition}

\newcommand{\R}{\mathbb R}
\newcommand{\C}{\mathbb C}
\newcommand{\K}{\cal K}
\renewcommand{\H}{{\mathbb H}}

\def\ra{\rightarrow}

\def\e{\emph}
\def\i{\infty}
\def\p{\partial}
\def\b{\begin}

\begin{document}

\title{    \flushleft{\bf{Quasiisometries between negatively curved Hadamard manifolds}}      }
\date{  }
\maketitle

\vspace{-5mm}

\noindent
Xiangdong Xie\newline
   Department of Mathematics,
   Virginia Tech,
   Blacksburg, VA 24060-0123\newline
   Email: xiexg@math.vt.edu



\pagestyle{myheadings}

\markboth{{\upshape Xiangdong Xie}}{{\upshape  Quasiisometries
between Hadamard manifolds}}


\vspace{3mm}

\noindent {\small {\bf Abstract.}
  Let $H_1$, $H_2$ be the universal covers of two compact
    Riemannian manifolds
  (of dimension $\not=4$)
  with negative sectional curvature. Then every quasiisometry between
   them lies at a finite distance from a bilipschitz homeomorphism.
     As a consequence,  every self quasiconformal map of
       a  Heisenberg group
       (equipped with the Carnot metric
    and  viewed as the ideal boundary of complex hyperbolic space)
       of dimension $\ge 5$
       extends to a  self quasiconformal map of the complex hyperbolic space.}



\vspace{3mm} \noindent {\small {\bf{Mathematics Subject
Classification (2000).}}  53C20,  20F65, 30C65.
}



\vspace{3mm} \noindent {\small {\bf{Key words.}}}   Quasisymmetric
extension, Hadamard manifolds,   quasiisometry,      bilipschitz.


\setcounter{section}{1}
\setcounter{subsection}{0}

\subsection{Introduction}\label{s1}

A classical result in quasiconformal analysis says that
 every quasisymmetric map from  $\R^n$
 to itself extends to a quasiconformal  map of  the upper half
  space  $$\R^{n+1}_+=\{(x_1, \cdots,  x_{n+1}):  x_{n+1}>0\}.$$
     Here   $\R^n$ is identified with
   the subset $\{(x_1, \cdots, x_{n+1}):  x_{n+1}=0\}$  of
   $\R^{n+1}$.
This result was first proved for $n=1$ by
  L. Ahlfors and A.  Beurling \cite{BA},  then
     for $n=2$ by L. Ahlfors  \cite{A}, for $n=3$
       by L.  Carleson   \cite{C}, and  finally   for all $n$ by
         Tukia and Vaisala   \cite{TV}.

 Recall  that   $\R^{n+1}_+$ is  the upper half space model
   for the real hyperbolic space and $\dot\R^n=\R^n\cup  \{\i\}$ is its ideal
   boundary.
The quasiconformal
     extension
  to  $\R^{n+1}_+$  turns out to be a bilipschitz homeomorphism in
  the hyperbolic metric.
     On the other hand,
       every self  quasiisometry  of   the hyperbolic space  $\R^{n+1}_+$
  induces a  self quasisymmetric map of  $\dot\R^n$ (equipped with the spherical metric), and
  conversely     every  self quasisymmetric map  of  $\dot\R^n$
   is  induced by a  self
   quasiisometry    of  $\R^{n+1}_+$ (see \cite{BS} for more  general statements).
  Hence what Tukia-Vaisala and others have proved is this:
    every  self  quasiisometry  of the hyperbolic space
$\R^{n+1}_+$   lies at a finite distance from a
 bilipschitz   homeomorphism.
   The main result of this paper generalizes this statement to
   negatively curved Hadamard manifolds.

    A Hadamard manifold
  is a simply connected complete Riemannian manifold with
    nonpositive sectional curvature.

\b{Th}\label{main} {Let $H_1$, $H_2$ be the universal covers of two
compact  Riemannian manifolds
  (of dimension $\not=4$)
with negative sectional curvature. Then every quasiisometry between
   them lies at a finite distance from a bilipschitz homeomorphism.
}

\end{Th}

The restriction on dimension  is due to the following  facts:
    a
$4$-dimensional topological manifold  might have
    more than
one lipschitz structures  or  none \cite{DS}, while
   every  topological $n$-manifold with $n\not=4$  has a unique
  lipschitz structure \cite{S}.

The unit ball  $B^n_\C$
  in  $\C^n$ is a model for the complex
hyperbolic space,  which is a negatively curved Hadamard manifold.
  The unit sphere in $\C^n$ has a sub-Riemannian structure
   coming from the complex structure of $\C^n$.
In this case the natural metric on the unit sphere is the associated
Carnot metric. 
  Theorem \ref{main}  applied to the case
    $H_1=H_2=B^n_\C$  yields a generalization of
     the classical quasiconformal extension result 
    to the complex hyperbolic case:

    \b{Cor}\label{ch}
    {Every  self  quasisymmetric map   of  the unit sphere
       (equipped with the Carnot metric)
       in $\C^n$  ($n\not=2$)
         extends to
       a self  quasiconformal  map  of the
          complex hyperbolic space  $B^n_\C$.

    }

    \end{Cor}

There are two general questions related to the results in this
paper.  It seems that   solutions to both questions require the use
  of controlled topology (for example the theorem of Chapman-Ferry
    \cite{CF}
  on small homotopy equivalences   or  its variants).
    The first general question
asks when a power quasisymmetric
 map between the boundaries of two Euclidean domains extend to a
 quasisymmetric map between the domains.  See Section \ref{s7} for a
 more precise formulation.

 The second general question asks
when   a quasiisometry lies at a finite distance from a bilipschitz
homeomorphism.
   Also see  Section \ref{s7} for a
 more precise formulation.
 For example, a natural question is whether the main theorem in this
 paper holds for all Hadamard manifolds, in particular, whether
 every self  quasiisometry  of $\R^n$ lies at a finite distance from
 a self bilipschitz homeomorphism of $\R^n$.
   The example of
  Dranishnikov-Ferry-Weinberger \cite{DFW} suggests   that
    the question is very subtle:
     they constructed two   
     uniformly contractible Riemannian manifolds
and a quasiisometry between them that is not at a finite distance
from any homeomorphism. The two manifolds they constructed  are
homeomorphic to the Euclidean space, but  the metrics are somehow
exotic.
  On the positive side,  Whyte \cite{W}   showed that
   every   quasiisometry between two
uniformly discrete non-amenable spaces of bounded geometry
  (for example, finitely
generated nonamenable groups with the word metric)
 lies    at a finite distance from a bilipschitz map.

Our proof follows the strategy of Tukia-Vaisala    \cite{TV}.
  In particular, we use the  boundary map induced by  the
  quasiisometry. This  means  that  our proof
   can not be generalized to nonpositively  curved spaces like
   $\R^n$: a  quasiisometry between CAT(0) spaces
  in general does not induce a boundary map.

 In Section \ref{s2} we review some basics about various
     maps
 and negatively curved  spaces. In Section  \ref{s3}
    we  replace the quasiisometry with a homeomorphism
      $F$, constructed using the boundary map of the quasiisometry.
       In general  the map $F$
  is not bilipschitz, but
       has very good  compactness property: both $F$ and $F^{-1}$
       are uniformly continuous. This is established in Sections
       \ref{s4} and \ref{s5}.  In Section \ref{s6} we
         modify $F$ to
        obtain a bilipschitz map. The arguments in Section \ref{s6}
           are  essentially
        due to Tukia-Vaisala   \cite{TV}.
  In the last Section we formulate some open questions.


\noindent {\bf{Acknowledgment}}. The author
  particularly thanks the
  Department of Mathematics at Virginia Tech for its generous
  support: the teaching load of one class per year is really a gift.
   The author would also
      like to thank Bruce Kleiner for
drawing his attention to the  paper of Block-Weinberger   \cite{BW}.

\subsection{Preliminaries}\label{s2}

\noindent
  {\bf{Various maps}}

A   bijection between two metric spaces
 $f: X\ra Y$ is  \e{$L$-bilipschitz} ($L\ge 1$)  if
  $$d(x,y)/L\le d(f(x), f(y))\le L\, d(x,y)$$
      for all
   $x,y\in X$.
     An embedding of metric spaces $f: X \ra Y$
      is \e{locally bilipschitz},  
         if each point $x\in X$
       has a neighborhood $U$ such that
         $f|_U$  is bilipschitz;  we say $f$ is \e{locally
         $L$-bilipschitz} if
           $f|_U$  is  $L$-bilipschitz.
   Let $L\ge 1$ and $A\ge 0$ be constants.
    A (not necessarily continuous) map $f:X \ra Y$
is  an  \e{$(L,A)$-quasiisometry} if the following holds:\newline
 (1)
$d(x,y)/L-A  \le d(f(x), f(y))\le L\, d(x,y)+A  $ for all
   $x,y\in X$;\newline
    (2)    For any $z\in Y$, there is some $x\in X$ with
     $d(z, f(x))\le A$.\newline
        If  $f:  I\ra Y$ is a map  defined on
         an interval  $I\subset \R$
  and satisfies condition (1) above, then we say $f$ is
    an \e{$(L, A)$-quasigeodesic}.

It is clear that a bilipschitz map is a quasiisometry, but the
converse is not true.

Let $\eta: [0,\i)\ra [0,\i)$ be a homeomorphism.
 A homeomorphism between metric spaces
$f: X\ra Y$ is \e{$\eta$-quasisymmetric} if for all
 distinct  triples   $x, y, z\in X$, we have
$$\frac{d(f(x), f(y))}{d(f(x), f(z))}\le \eta\left(\frac{d(x,y)}{d(x,z)}\right).$$
A homeomorphism $f: X\ra Y$ is quasisymmetric  if it is
$\eta$-quasisymmetric  for some $\eta$.
 A quasisymmetric map is called a \e{power quasisymmetric map} if it
   is $\eta$-quasisymmetric   and  $\eta$ has the form
    $\eta(t)=C t^\alpha$, where $C>0$ and $\alpha>0$ are constants.
 Recall that a quasisymmetry between connected metric spaces is a
 power   quasisymmetry (see Theorem 6.14 in \cite{V}).


\vspace{3mm}

\noindent
{\bf{Gromov hyperbolic spaces}}

Let $X$ be a geodesic metric space.
  We assume $X$ is \e{proper}, that is, all closed balls in $X$ are
  compact.
    Let $\delta\ge 0$. We say $X$ is \e{$\delta$-hyperbolic},  if for
     any $x, y, z\in X$, and any geodesics $xy$, $yz$, $zx$ between
     them,  $yz$ is contained in the $\delta$-neighborhood of
      $xy\cup xz$.
 A  metric space is \e{Gromov hyperbolic} if it is $\delta$-hyperbolic
 for some $\delta$.

  A Gromov hyperbolic geodesic space
$X$ has an ideal boundary $\p X$:  by definition, $\p X$ is the set
 of equivalence classes of geodesic rays in $X$, where two rays are
 equivalent if the Hausdorff distance between them is finite.
   There is a natural topology on $\overline X:=X\cup \p X$,  
    in which  $\overline X$ is compact and
  $X$ is an open dense subset of $\overline X$.
    Let  $X$ and $Y$ be Gromov hyperbolic spaces.
        Then  every  quasiisometry
       $f: X\ra Y$  induces a boundary
     map  $\p f: \p X\ra \p Y$, which     is a
     homeomorphism.
       Moreover,  $\p f$ is a power quasisymmetry with respect to
       the so-called visual metrics on the Gromov boundary;
         conversely,
       if $X$ and $Y$ satisfy some mild conditions, then every
         power quasisymmetry $\p X\ra \p Y$ is induced
           by a quasiisometry. See \cite{BS} for more details. See
           also Lemma \ref{l3.1}  for a  statement in the
           Hadamard manifold case.

Let  $X$ be a $\delta$-hyperbolic geodesic space,
  $x, y, z\in X$, and     $xy$, $yz$, $zx$
      geodesics    between them.
    We say  $xy\cup yz\cup zx$  is a triangle.
    Let $C\ge 0$.
We say
 a point $w\in X$ is  a   \e{$C$-quasicenter} of  $xy\cup yz\cup zx$
  if  $d(w, xy), d(w, xz), d(w, yz)\le C$.
     There is a constant $C'=C'(\delta, C)$ with the following
     property:  for any triangle
     in $X$ and   any its  two $C$-quasicenters
      $w_1$, $w_2$,  
         the inequality
        $d(w_1, w_2)\le C'$ holds.

Quasigeodesics in Gromov hyperbolic spaces have the so-called
stability property.  For any $\delta\ge 0$, any $L\ge 1$, $A\ge 0$,
there is a constant $C=C(L, A, \delta)$ with the following
  property:  for any $\delta$-hyperbolic space  $X$,
    any two $(L, A)$-quasigeodesics   $\gamma_1: I_1\ra X$,
      $\gamma_2: I_2\ra X$
  with the same endpoints, the inequality
    $HD(\gamma_1(I_1), \gamma_2(I_2))\le C$ holds.
Here we use the  notation $HD_d(A, B)$ to denote the Hausdorff
distance between two subsets $A, B\subset X$ of a metric space $(X,
d)$;  we often write $HD(A, B)$ if the metric in question is clear.

\vspace{3mm}

\noindent
  {\bf{Hyperbolic trigonometry}}


  For $\lambda>0$,  let
   $\H^2(-\lambda^2)$ be  the hyperbolic plane with constant sectional
curvature $-\lambda^2$.  We  abbreviate $\H^2(-1)$ by $\H^2$. 

Let $\Delta$  be a triangle in $\H^2(-\lambda^2)$.  Denote the three
angles by $A, B, C$ and
  the lengths  of
 their opposite   sides by $a, b, c$.
If  the angle $C$  is a right angle,  then (\cite{G}, p.24)
$$\cosh (\lambda a)=\frac{\cos A}{\sin B}.$$

\vspace{3mm}

\noindent
  {\bf{Hadamard manifolds}}

    Let $H$ be a Hadamard manifold.
  A  classical theorem of Hadamard says that for every $x\in H$,
    the exponential map $\exp_x:  T_x H\ra H$  is a diffeomorphism
      from the tangent space $T_x H$ onto $H$.
   The  ideal boundary $\p H$  of $H$
    is   defined in the same way as
   for Gromov hyperbolic spaces.
   There is the so-called cone  topology on $\overline H:=H\cup \p H$,
       in which   $\p H$ is homeomorphic to a sphere
   and
    $\overline H$ is homeomorphic to a closed
   ball.  

The distance function on  $H$ is convex: for any two geodesics
  $c_1, c_2: \R\ra H$,  the function $f(t):=d(c_1(t), c_2(t))$ is
  convex.   In particular,  if $c_1(\i)=c_2(\i)$,  then
$d(c_1(t), c_2(t))$ is  decreasing.

For each $x\in H$, let $S_x\subset T_x H$ be the unit tangent sphere
of $H$ at $x$.  There is a map $L_x: S_x\ra \p H$,
  where  for each $v\in S_x$,
  $L_x(v)$  is  the equivalence class containing
      the ray  starting at $x$ with initial
  direction    $v$.  The map $L_x$ is a homeomorphism.  It follows
  that for each equator (intersection of $S_x$
    with    a hyperplane of $T_x H$)
      in $S_x$, its image under $L_x$
  is a codimension 1 sphere separating $\p H$ into two balls.

For $x\in H$ and $y, z\in \overline H\backslash\{x\}$,
  the angle $\angle_x(y,z)$  is defined to be the angle
    in the tangent space $T_x H$ between the
  initial directions of $xy$ and $xz$.
   This is a continuous function on $y, z\in \overline
   H\backslash\{x\}$:
     if $\{y_i\}_{i=1}^\i,  \{z_i\}_{i=1}^\i\subset
\overline
   H\backslash\{x\}$  are two sequences with
    $y_i\ra y$ and $z_i\ra z$, then
     $\angle_x(y_i, z_i)\ra \angle_x(y,z)$.

{\bf{From now on, we shall  assume  that
  the   sectional curvature   of  $H$  satisfies
    $-\lambda^2\le K\le -1$, where $\lambda\ge 1$ is some
      fixed   constant.}}
Then $H$ is a $\delta_0$-hyperbolic space, where $\delta_0=\log 3$.
  See \cite{CDP}, p.12.

  There is a family of visual metrics  on  $\p H$:
     there exists  a universal constant $C_0$, such that  for each $x\in H$, there is  a metric
     $d_x$ on $\partial H$ satisfying
       $$\frac{1}{C_0} e^{-d(x, \xi\eta)}\le d_x(\xi, \eta)\le C_0 e^{-d(x,
       \xi\eta)}$$
           for all $\xi\not=\eta\in \partial H$. See \cite{B}
       Section 2.5.  Here $\xi\eta$ is the geodesic connecting $\xi$
       and $\eta$.

Let $\xi_1\not= \xi_2\in \p H$  and set
  $[\xi_1, \xi_2]=\xi_1\xi_2\cup \{\xi_1, \xi_2\}$.
    The orthogonal projection
$P_{\xi_1\xi_2}: \overline H\ra [\xi_1, \xi_2]$  is defined as
follows: for $x\notin [\xi_1, \xi_2]$, $P_{\xi_1\xi_2}(x)$ is the
unique point $w$ on $\xi_1\xi_2$ such that $wx$ is perpendicular to
$\xi_1\xi_2$,   and for $x\in [\xi_1, \xi_2]$,
 $P_{\xi_1\xi_2}(x)=x$.

  For $x, y, z\in H$,  let
        $\Delta(x, y, z)$  be the  triangle
  with vertices $x, y, z$.
         A triangle $\Delta(x_\lambda, y_\lambda, z_\lambda)$ in $\H^2(-\lambda^2)$ is a
         comparison triangle of $\Delta(x,y,z)$
           if $d(x,y)=d(x_\lambda, y_\lambda)$,
            $d(y, z)=d(y_\lambda, z_\lambda)$
              and $d(z, x)=d(z_\lambda, x_\lambda)$.
  Comparison triangle always exists and is  unique up to isometry.
    The comparison angle
  $\tilde{\angle_x}(y,z)$  is defined to be
  $\angle_{x_\lambda}(y_\lambda , z_\lambda)$.
 For  $p\in xy$ and $q\in xz$,  points $p_\lambda\in x_\lambda y_\lambda$
  and $q_\lambda\in x_\lambda z_\lambda$ are said to correspond to $p$ and $q$
    if $d(p,x)=d(p_\lambda, x_\lambda)$ and $d(q, x)=d(q_\lambda, x_\lambda)$.

  Let  $\Delta(x, y, z)$  be a triangle in $H$.
   Let $\Delta(x_\lambda, y_\lambda, z_\lambda)$ and $\Delta(x_1, y_1, z_1)$ respectively
     be the  comparison triangles of $\Delta(x, y, z)$ in
     $\H^2(-\lambda^2)$ and $\H^2$.
       Then we have the following:\newline
       (1) $\angle_{x_\lambda}(y_\lambda, z_\lambda)\le \angle_x(y,z)\le \angle_{x_1}(y_1,
       z_1)$;  \newline
  (2) For any $p\in xy$, $q\in xz$, if $p_\lambda\in x_\lambda y_\lambda$,  $q_\lambda\in
  x_\lambda z_\lambda$  and $p_1\in x_1y_1$, $q_1\in x_1z_1$  correspond to
   $p$, $q$,  then
     $d(p_\lambda, q_\lambda)\le d(p,q)\le d(p_1, q_1)$.\newline
       See \cite{BH}  p.161, p.169-173
  and \cite{CE} p.42  for more details.

\b{Le}\label{l2.1} {Given any $\epsilon>0$, there is a constant
 $\epsilon_1'=\epsilon_1'(\epsilon)>0$  with the following property:
for any three points $x,y,z\in \H^2$, if
  $d(x,y)\ge\epsilon$  and $\angle_x(y,z), \angle_y(x,z)\ge \pi/4$,
    then  $\angle_x(y,z)+ \angle_y(x,z)\le \pi-\epsilon_1'$.

}

\end{Le}

\b{proof} Let $m$ be the  midpoint of $xy$ and
  $x'\in mx$ and $y'\in my$  with $d(m, x')=d(m, y')=\epsilon/2$.
 Let $S$ be the unique  circle   in $\H^2$  satisfying the
following: \newline (1)
  the point $z$ and the center  $p$ of $S$ lie at  the  same side of
  $xy$; \newline
   (2)
  $S$  is  tangent to $xy$ at the midpoint  $m$ of $xy$; \newline
 (3)
  $\angle_{x'}(p, y)=\pi/8$.  \newline
    Then the ball $B$ inside $S$ is
 contained in the triangle $\Delta(x,y,z)$, due to our assumption on
 the angles
   $\angle_x(y,z)$, $\angle_y(x,z)$.  Hence the area   $A$ of
$\Delta(x,y,z)$ is at least the area of $B$,  which is a constant
depending only on $\epsilon$.  Now the lemma follows from the
Gauss-Bonnett  formula:
  $\angle_x(y,z)+\angle_y(x,z)+\angle_z(x,y)=\pi-A$.

\end{proof}

\b{Le}\label{l2.2} {Let $c_i:  [0,\i)\ra H$ \e{($i=1, 2$)}
  be two   equivalent  rays. Set
   $x=c_1(0)$, $y=c_2(0)$  and $\xi=c_1(\i)$.   Suppose $d(x,y)\ge \epsilon$ and
     $\angle_x(y, \xi), \angle_y(x,\xi)\ge 3\pi/8$.
       Then $\angle_x(y, \xi)+ \angle_y(x,\xi)\le \pi-\epsilon_1'$,
         where $\epsilon_1'$  is the constant in Lemma \ref{l2.1}.

}

\end{Le}

\b{proof} For any $0<\eta<\pi/8$, pick a point $z\in y\xi$ with
 $|\angle_x(y, z)-\angle_x(y, \xi)|\le \eta$.
   Consider the comparison triangle of $\Delta(x,y,z)$
     in $\H^2$.  Since $H$ has sectional curvature $K\le -1$,
      we have $\tilde{\angle_x}(y,z)\ge \angle_x(y,z)
        \ge \angle_x(y, \xi)-\eta\ge \pi/4$ and
          $\tilde{\angle_y}(x, z)\ge \angle_y(x,z)\ge 3\pi/8$.
            It now follows from Lemma \ref{l2.1} and the assumption
             $d(x,y)\ge \epsilon$ that
               $\tilde{\angle_x}(y,z)+\tilde{\angle_y}(x, z)\le
               \pi-\epsilon_1'$.
               Hence
$\angle_x(y, \xi)+ \angle_y(x,\xi)\le \pi-\epsilon_1'+\eta$.
  The lemma follows by letting $\eta\ra 0$.

\end{proof}

\b{Le}\label{l2.3} {Given any $\epsilon>0$, $\epsilon_1'>0$,  there
is some
   $\epsilon_1''=\epsilon_1''(\epsilon, \epsilon'_1)>0$ with the following property:
     for any $x,y,z\in H$, if $d(x,z)\ge \epsilon$ and
     $d(x,y)\le \epsilon_1''$, then
      $\angle_z(x,y)\le \epsilon_1'/10$.

}

\end{Le}

\b{proof} The statement  clearly holds for $\H^2$. The statement for
$H$   follows  by comparison.

\end{proof}

\b{Le}\label{l2.4}
  {There is an absolute constant $C_1$ with the
 following property:
 for any $x\in H$ and $y, z\in \overline H\backslash \{x\}$,
     if $\angle_x(y, z)\ge \pi/2$, then $d(x, yz)\le C_1$;
        furthermore,  $HD(yz, xy\cup xz)\le C_1$.

  }

  \end{Le}

  \b{proof}
It suffices to prove the first claim  for $y, z\in H$.
   It  clearly holds for $\H^2$, and
    the general case follows by comparison.
The second claim then follows from the
 convexity of distance function.

  \end{proof}






The following lemma says that for $x\in H$ and $\xi, \eta\in \p H$,
  $d_x(\xi, \eta)$ is small if
      and only  if  $\angle_x(\xi,\eta)$ is small.

\b{Le}\label{l2.5} {For any $\epsilon>0$, there is some
$\delta=\delta(\epsilon)>0$
  with the following properties:  for
 $x\in H$ and $\xi, \eta\in \p
H$,  \newline
  (1)  
    if  $d_x(\xi,\eta)<\delta$, then  $\angle_x(\xi,\eta)<
    \epsilon$;\newline
      (2) 
    if  $\angle_x(\xi,\eta)<
    \delta$,  then
    $d_x(\xi,\eta)<\epsilon$.  

}

\end{Le}

\b{proof} 
  Let $z\in \xi\eta$ be the point
such that $xz$ is perpendicular to $\xi\eta$.

  (1)   Fix $\epsilon>0$.
 By Lemma \ref{l2.4},  $d(z, x\xi), d(z, x\eta)\le C_1$.
  There is some $b=b(\epsilon)>0$ such that
    $\angle_x(z, \xi)<\epsilon/2$   whenever $d(x, z)\ge b$.
      Similarly for $\angle_x(z, \eta)$.
        Now the claim follows from
         $d_x(\xi, \eta)\ge e^{-d(x, \xi\eta)}/{C_0}=e^{-d(x,
         z)}/{C_0}$.


(2)   
 For $n\ge 1$,
   let $\xi_n\in z\xi$,  $\eta_n\in z\eta$  be points at distance
   $n$ from  $z$. Then  $\angle_x(\xi_n, \eta_n)<2\angle_x(\xi,
   \eta)$  for sufficiently large $n$.
Let $\Delta(\bar x, \bar {\xi}_n, \bar{\eta}_n)$  be a
  comparison triangle of $\Delta(x, \xi_n, \eta_n)$  in
  $\H^2(-\lambda^2)$  and $\bar z_n\in \bar {\xi}_n\bar{\eta}_n$ the
  point corresponding to $z$.
    Since
      $H$ has curvature $K\ge -\lambda^2$,  we have
          $\angle_x(\xi_n, \eta_n)\ge
    \angle_{\bar x}(\bar {\xi}_n, \bar{\eta}_n)$  and
  $d(x, z)\ge d(\bar x, \bar{z}_n)$.
  Notice   that
    $d(\bar x, \bar{z}_n)$ is bounded  above by a number
  independent of $n$ and
  $d(\bar x, \bar{\eta_n})\ra \i$ as $n\ra \i$.  It follows that
  $\angle_{\bar{\eta}_n}(\bar x, \bar{\xi}_n)\ra 0$.
  Similarly,  $\angle_{\bar{\xi}_n}(\bar x, \bar{\eta}_n)\ra 0$.
   Hence  for sufficiently large $n$,
    the projection   $w$ of $\bar x$ on
   $\bar{\xi}_n\bar{\eta}_n$  lies in the interior of
$\bar{\xi}_n\bar{\eta}_n$,  $\bar x w$ is perpendicular to
$\bar{\xi}_n\bar{\eta}_n$  and  $\angle_{\bar x}(w, \bar{\eta}_n)\le
\angle_{\bar x}(\bar {\xi}_n, \bar{\eta}_n)<2\angle_x(\xi,
   \eta)$. 
    Also notice $d(\bar x, w)\le d(\bar x,
    \bar z_n)\le d(x,z)=d(x, \xi\eta)$.

  We  may assume $\angle_x(\xi,
   \eta)\le \pi/4$.
      Now we  have
       $$\cosh(\lambda\, d(x, \xi\eta))\ge \cosh(\lambda\,d(\bar x,
   w))=\frac{\cos\angle_{\bar{\eta}_n}(\bar x,
       w)}{\sin \angle_{\bar x}(w, \bar{\eta}_n)}\ge
\frac{\cos\angle_{\bar{\eta}_n}(\bar x,
       \bar{\xi}_n)}{\sin 2\angle_x(\xi,
   \eta)}\ra \frac{1}{\sin 2\angle_x(\xi,
   \eta)}.$$
       Hence
$\cosh(\lambda\, d(x, \xi\eta))\ge \frac{1}{\sin 2\angle_x(\xi,
   \eta)}$.
  Now the claim follows since
  $d_x(\xi,  \eta)\le C_0e^{-d(x, \xi\eta)}$.


\end{proof}

\b{Le}\label{l2.6} {Given any $\epsilon>0$, there is some
$\delta=\delta(\epsilon, \lambda)>0$ with the following property:
    for any  three distinct points $p, \xi,\eta\in \p H$,  any $x\in
  p\xi$,   \newline
    (1)  if $d_x(\xi,\eta)<\delta$, then
   $d(x, p\eta)<\epsilon$;  \newline
     (2)  if  $d(x,p\eta)<\delta$, then $d_x(\xi, \eta)<\epsilon$.

}

\end{Le}

\b{proof}
  (1)  Fix $\epsilon>0$.  Let $\epsilon'_1$ be the constant in
  Lemma \ref{l2.1}.  By Lemma \ref{l2.5}, there is a constant
  $\delta=\delta({\epsilon'_1}/10)>0$ such that if $d_x(\xi,\eta)<\delta$, then
  $\angle_x(\xi,\eta)<\epsilon'_1/10$.
    Then $\angle_x(p, \eta)\ge \pi-\epsilon'_1/10$.
  Let $z$ be the projection of $x$ on $p\eta$.
  Then $\angle_x(z, \eta), \angle_x(z, p)< \pi/2$
and   $\angle_x(z, \eta)+ \angle_x(z, p)\ge \angle_x(p, \eta)\ge
\pi-\epsilon'_1/10$.
  It follows that
$\angle_x(z, \eta), \angle_x(z, p)\ge \pi/2-\epsilon'_1/10 $.
  Now Lemma \ref{l2.2}  applied to $xp$ and $zp$
  implies that $d(x, p\eta)=d(x,z)< \epsilon$.

(2) {\bf{Claim:}}  for any $\epsilon>0$, there is some
$\delta=\delta(\epsilon, \lambda)>0$ such that
  for  $x\in H$,  $\xi,\eta\in \p H$,  if $d(x, \xi\eta)<\delta$,
  then $\angle_x(\xi, \eta)> \pi-\epsilon$.
    To see this, first notice that for any
 $\epsilon>0$, there is some $\delta=\delta(\epsilon, \lambda)>0$ with the following property:
   for any $x, y,z\in \H^2(-\lambda^2)$,  if $d(x,y)=d(x,z)\ge 10$  and
     $d(x, yz)\le \delta$, then $\angle_x(y,z)> \pi-\epsilon$.
  Now for  $x\in H$,  $\xi,\eta\in \p H$  with  $d(x, \xi\eta)<\delta$,
 choose $y\in x\xi$ and $z\in x\eta$   such that
  $d(x,y)=d(x,z)\ge 10$  and
$d(x, yz)<\delta$.  Let $\Delta(x_\lambda,y_\lambda,z_\lambda)$ be
  a
comparison triangle of $\Delta(x,y,z)$ in $\H^2(-\lambda^2)$.
  Then $d(x_\lambda, y_\lambda z_\lambda)\le d(x,yz)$ and $\angle_x(y,z)\ge
  \angle_{x_\lambda}(y_\lambda, z_\lambda)$.
    Since $\angle_x(\xi,\eta)=\angle_x(y,z)$, the claim follows.

Now fix $\epsilon>0$.
 By Lemma \ref{l2.5}, there is some $\epsilon''=\epsilon''(\epsilon)>0$ such that if
  $\angle_x(\xi,\eta)<\epsilon''$, then $d_x(\xi,\eta)<\epsilon$.
   By the above claim, there is some $\delta=\delta(\epsilon'',\lambda)>0$ such that
     if   $d(x, p\eta)<\delta$, then $\angle_x(p,
     \eta)>\pi-\epsilon''$.   Now
     assume  $d(x, p\eta)<\delta$.  Then
      $\angle_x(\xi,\eta)=\pi-\angle_x(p, \eta)<\epsilon''$
  and   hence  $d_x(\xi,\eta)<\epsilon$.

\end{proof}

\b{Le}\label{l2.7} {
  Given any $\epsilon>0$, there is some $\delta=\delta(\epsilon, \lambda)>0$ with the following
  property:
    for  $x,y\in H$ and  $\xi\in \p H$, if
       $\angle_x(y, \xi)<\pi/2$, $\angle_y(x, \xi)\le
      \pi/2$  and   $d(x,y)\le \delta$,   then
         $\angle_x(y, \xi)>\pi/2-\epsilon$.

}

\end{Le}

\b{proof} Let $z_n\in y\xi$ be the point at distance  $n$ from $y$.
  Then $\angle_x(y, \xi)=\lim_{n\ra \i}\angle_x(y, z_n)$.
 Consider a comparison triangle  $\Delta(x', y', z'_n)$
  of $\Delta(x, y,z_n)$  in $\H^2(-\lambda^2)$.
  Then
  $$\angle_{x'}(y', z'_n)\le \angle_x(y, z_n)\le \pi/2 \;\;\; \text{and}  \;\;
    \angle_{y'}(x', z'_n)\le  \angle_y(x, z_n)= \angle_y(x, \xi)\le \pi/2.$$
  Let $w'$ be the projection of $z'_n$ on $x'y'$.
    Now  by hyperbolic trigonometry
      $$\sin \angle_{x'}(y', z'_n)=\frac{\cos
      \angle_{z'_n}(w', x')}{\cosh (\lambda\,d(w', x'))}\ge
\frac{\cos
      \angle_{z'_n}(w', x')}{\cosh (\lambda\,d(x,y))}.$$
        Since $ \angle_{z'_n}(w', x')\ra 0$ as $n\ra \i$,  
    we  have
  $$\sin \angle_x(y, \xi)=\lim_{n\ra \i} \sin \angle_x(y, z_n)\ge \liminf_{n\ra \i}\sin\angle_{x'}(y', z'_n)\ge
\frac{1}{\cosh (\lambda\,d(x,y)))}.$$
   It follows that there is a function $g(t)$ satisfying
     $g(t)\ra \pi/2$ as $t\ra 0$  such that
       $\angle_x(y, \xi)\ge g(d(x,y))$.  The lemma follows from
       this.

\end{proof}

\subsection{Constructing   the map  $F$}\label{s3}

Let $H_1$, $H_2$ be the universal covers of two compact Riemannian
manifolds   with negative sectional curvature,
  and
    $f: H_1\ra
H_2$ a  quasiisometry.
 In this section we construct a map
  $F: H_1\ra H_2$ which has the same boundary map
    as $f$.    In general $F$ is not a bilipschitz
    homeomorphism  and we shall modify $F$ in Section \ref{s6}
      to  obtain a bilipschitz map.

    Notice that, after rescaling  the metrics on $H_1$ and $H_2$, we
    may assume   that  their sectional curvature
        satisfies
$-\lambda^2\le K\le -1$  for some constant $\lambda\ge 1$.
  We shall assume this from now on.

The  Hadamard manifolds
  $H_1$, $H_2$  are Gromov hyperbolic spaces,
  and
$f$ induces
 a  boundary map $\partial f: \partial H_1\ra \partial H_2$.
       We set
   $\xi':=\partial f(\xi)$ for any
     $\xi\in \partial H_1$.
        Fix a point  $p\in \partial H_1$.
       For each $\xi\in \partial H_1\backslash\{p\}$,
     the map $F$ shall send  the geodesic line $p\xi\subset H_1$ into the
     geodesic line $p'\xi'\subset H_2$.
 The  map  $F$  depends on the  point $p$.  However, we shall
 suppress this information to simplify the notation.

 Let $\xi\in\partial H_1\backslash\{p\}$  and
  $x\in p\xi$.
Let $F_x\subset S_x H_1$  be the set of unit tangent vectors at $x$
that are perpendicular to  $p\xi$.  Then $F_x$ is an equator in $S_x
H_1$.
  Set $E_x=L_x(F_x)\subset \partial H_1$, where $L_x: S_x H_1\ra \partial H_1$
     is the  homeomorphism   defined in Section \ref{s2}.
      Then  $E_x$ is a codimension 1 sphere in $\partial H_1$
        separating  $\partial H_1$ into two open balls
that contain  $p$  and  $\xi$  respectively.
        We let $F(x)$ be the
point in
  $P_{p'\xi'}(\partial f(E_x))$ that is closest to $p'$;
    that is,   if $c: \R\ra H_2$ is an arc-length parametrization of
      $\xi'p'$ from $\xi'$ to $p'$  and
      $$t_x:=\sup\{t\in \R:
         c(t)=P_{p'\xi'}(\partial f(\beta)) \;\; {\text{for some}}
         \;\;   \beta\in E_x\},$$
              then
        $F(x)=c(t_x)$.
 The compactness of $\partial f(E_x)$  and the continuity of the projection
  $ P_{p'\xi'}$  imply that
         $F(x)\in  P_{p'\xi'}(\partial
        f(E_x))$.

Notice that $\partial f(E_x)$  separates $\partial H_2$ into two
open balls  that contain  $p'$  and  $\xi'$  respectively.
    This property implies that $F$ is injective
  along $p\xi$. Since $\partial f$ is a homeomorphism,
    $F$ is injective.  In the next two sections we shall show that
    $F$ has good compactness property:  both $F$ and $F^{-1}$ are
      uniformly continuous.
 Below we  first prove some preliminary results.

Let $f:H_1\ra H_2$  be  an   $(L, A)$-quasiisometry.
 The point of the following
   lemma  is that the control function $\eta$
   is independent of the base point  $x$.

\b{Le}\label{l3.1} {Let $p,q,r\in \p H_1$ be three distinct points.
  Set $x=P_{pq}(r)$ and $x'=P_{p'q'}(r')$.
 Then
  $\partial f: (\p H_1, d_x)\ra (\p H_2, d_{x'})$  is an
  $\eta$-quasisymmetry, where $\eta(t)=C t^{\frac{1}{L}}$
    and $C$ depends only on $L$ and $A$.

}

\end{Le}

\b{proof} Notice that $xr$ is perpendicular to $pq$.  By Lemma
\ref{l2.4}
  there is an absolute constant  $C_1$   with
  $d(x, pr), d(x, qr)\le  C_1$.
       Similarly  $d(x', p'r'), d(x', q'r')\le  C_1$.
       It follows that $d(f(x),  f(pr)),  d(f(x), f(qr))\le LC_1+A$.
Since $f$ is a $(L, A)$-quasiisometry, the images of $pr$ and $qr$
under $f$  are  $(L,A)$-quasigeodesics.  By the stability  of
quasigeodesics, 
    there is a constant
$C_2=C_2(L, A)$,  such that $HD(p'r', f(pr))\le C_2$ and
 $HD(q'r',  f(qr))\le C_2$.
  Hence $d(f(x), p'r'), d(f(x), q'r'), d(f(x), p'q')\le C_2+LC_1+A$.
 That is, $f(x)$ is a   $C_3$-quasicenter  of the three points
   $p', q',  r'$, where $C_3=C_2+LC_1+A$.  Since $x'$ is also a  $C_3$-quasicenter of   $p',q',  r'$,
     there is a constant $C_4=C_4(L,A)  $   
         such
     that
       $d(x', f(x))\le C_4$.

Now for any $\xi,\eta\in \p H_1$, we have
   \b{align*}
   d(x',
\xi'\eta')\ge  &  d(f(x), \xi'\eta')-d(x', f(x))\ge d(f(x),
f(\xi\eta))-HD(\xi'\eta', f(\xi\eta))-C_4\\
     \ge   &  d(f(x), f(\xi\eta))-C_2-C_4\ge d(x, \xi\eta)/L-A-C_2-C_4.
\end{align*}
 It
follows that
  $$d_{x'}(\xi',   \eta')\le C_0 e^{-d(x', \xi'\eta')}\le
C_0e^{A+C_2+C_4}(e^{-d(x, \xi\eta)})^{\frac{1}{L}} \le C (d_x(\xi,
\eta))^{\frac{1}{L}},$$
  where $C=C_0^{1+ {\frac{1}{L}}}
e^{A+C_2+C_4}$ depends only on $L$ and $A$.

\end{proof}

For any $x\in H_1$,  we let $q_x\in \p H_1$ be the unique point such
that  $x\in pq_x$.


\b{Le}\label{l3.2} {There are two absolute constants $B_2\ge B_1>0$
such that
  $B_1\le  d_x(\xi, q_x)\le B_2$
  for all $x\in H_1$  and all $\xi\in E_x$.

}

\end{Le}

\b{proof} Let  $x\in H_1$ and $\xi\in E_x$.
 Then $x\xi$ is perpendicular to $pq_x$.  By Lemma \ref{l2.4}
    $d(x, \xi q_x)\le C_1$.
 The lemma  now  follows from    
  $e^{-d(x, \xi_1\xi_2)}/{C_0}\le d_x(\xi_1, \xi_2)\le C_0 e^{-d(x,
  \xi_1\xi_2)}$.

\end{proof}

From now on we set $x':=F(x)$ for $x\in H_1$.

\b{Le}\label{l3.3} {There are two constants $B_4\ge B_3>0$ that
depend only on $L$ and $A$ such that
   $B_3\le d_{x'}(q_x', \beta')\le B_4$
  for all $x\in H_1$ and all $\beta\in E_x$.

}

\end{Le}

\b{proof}
  There is some $\xi\in E_x$ such that
    $\xi'\in E_{x'}$. By Lemma \ref{l3.2},  we have
     $B_1\le d_{x'}(q'_x, \xi')\le B_2$.
 On the other hand, Lemma \ref{l3.1}  says that
 $\partial f: (\p H_1, d_x)\ra (\p H_2, d_{x'})$  is an
  $\eta$-quasisymmetry.   For any $\beta\in E_x$,
   we have
   $$\frac{d_{x'}(q'_x, \beta')}{d_{x'}(q'_x, \xi')}\le
   \eta\left(\frac{d_x(q_x, \beta)}{d_x(q_x, \xi)}\right)\le
   \eta\left(\frac{B_2}{B_1}\right).$$
  Hence $d_{x'}(q'_x, \beta')\le \eta(\frac{B_2}{B_1})d_{x'}(q'_x,
  \xi')\le B_2\eta(\frac{B_2}{B_1})$.
   Similarly by considering
     $\frac{d_{x'}(q'_x, \xi')}{d_{x'}(q'_x, \beta')}$  we obtain
$d_{x'}(q'_x, \beta')\ge   \frac{B_1}{\eta(\frac{B_2}{B_1})}$.

\end{proof}

\b{Le}\label{l3.4} { There is a constant
 $B_5=B_5(L,A)$    with the following
 property:
    for any $x\in H_1$, the projection
      $P_{p'q'_x}(\p f(E_x))$  of $\p f(E_x)$ on
        the geodesic $p'q'_x$  is a closed segment with length
          at most $B_5$.

}

\end{Le}

\b{proof} $P_{p'q'_x}(\p f(E_x))$  is a closed segment since
  the projection is continuous and
  $\p f(E_x)$ is compact and  connected.
By the definition of $x'$ we know that   $x'$ is the endpoint of the
segment that is closer to $p'$.
    Let   $\eta'\in \p f(E_x)$. Denote by $z$ the projection of
  $\eta'$ on $p'q'_x$ and $w$ the projection of $x'$ on $q'_x\eta'$.
    Since
      $B_3\le d_{x'}(q'_x, \eta')\le C_0 e^{-d(x', q'_x\eta')}$,
          we  have    $d(x', w)\le C$  for some
         $C=C(L, A)$.  Lemma \ref{l2.4} applied to
          $z, p', \eta'$ implies that
           $d(x', p'\eta')\le C_1$.
            It follows that $w$ is a $(C+C_1)$-quasicenter of
             $p', \eta', q'_x$.
              Lemma \ref{l2.4} also implies  that
                $z$ is a
              $C_1$-quasicenter  of  $p', \eta', q'_x$.
Hence $d(z,w)\le C'$ for some $C'=C'(L, A)$.
 Now $d(x', z)\le d(x', w)+d(w, z)\le C+C'$.


\end{proof}

\subsection{$F$ is uniformly continuous}\label{s4}

In this Section we prove that $F$ is uniformly continuous.
  It follows  that $F$ is a homeomorphism.
    For any $\epsilon>0$, we need to find $\delta>0$ such that
   $d(F(x), F(y))< \epsilon$
     for all $x,y\in H_1$  satisfying  $d(x,y)< \delta$.
 This is achieved in Lemmas \ref{l4.1}--\ref{l4.4}.

Recall our notation: $q'= \p f(q)$ for $q\in \p H_1$ and $x'=F(x)$
for $x\in H_1$.
  For simplicity, for $\alpha, \beta>0$,
   the notation $\beta=\beta(\alpha)$ means that
      the number $\beta$
   depends on $\alpha$ and possibly  $L$, $A$, $\lambda$, but
   nothing else.

\b{Le}\label{l4.1} {Given any $\epsilon>0$, there exists
$\epsilon_1=\epsilon_1(\epsilon)>0$ with the following property:
  for any $x,y\in H_1$,  if  $d_{x'}(q'_x, q'_y)<\epsilon_1$  and
    $HD_{d_{x'}}(\p f(E_x), \p f(E_y))<\epsilon_1$,  then
    $d(x',y')<\epsilon$.

}

\end{Le}

\b{proof} Fix $\epsilon>0$.

 \noindent
  {\bf{Claim:}}
       there exists $\delta_1=\delta_1(\epsilon)>0$ with the
following property:
  for any $x,y\in H_1$,  if  $d_{x'}(q'_x, q'_y)<\delta_1$  and
    $HD_{d_{x'}}(\p f(E_x), \p f(E_y))<\delta_1$,  then
 for any  $z\in P_{p'q'_y}(\p f(E_{y}))$, there is some $w\in
 P_{p'q'_x}(\p f(E_{x}))$  such that  $d(z,w)<
 \epsilon/2$.

 Assuming the claim, we first finish the proof of the lemma.
  Set $\epsilon_1=\delta_1/({{C^2_0}e^{B_5+\epsilon}})$.  Let
$x,y\in H_1$  and   assume   $d_{x'}(q'_x, q'_y)<\epsilon_1$,
    $HD_{d_{x'}}(\p f(E_x), \p f(E_y))<\epsilon_1$.
  By Lemma \ref{l3.4},  $P_{p'q'_x}(\p f(E_{x}))$
      has length  at most
  $B_5$.    It follows from  the claim  that    $d(x', y')<  \epsilon/2+B_5$.
    The inequality
     $e^{-d(z, \xi\eta)}/{C_0}\le  d_z(\xi, \eta)\le C_0e^{-d(z, \xi\eta)}$
  now implies
 $$d_{y'}(q'_x, q'_y)<\delta_1 \;\; {\text{and}}\;\;
    HD_{d_{y'}}(\p f(E_x), \p f(E_y))<\delta_1.$$
       Now  switching  the role of $x$ and $y$ and applying
          the claim, we see that
 the Hausdorff distance between
 $P_{p'q'_x}(\p f(E_{x}))$ and
    $P_{p'q'_y}(\p f(E_{y}))$
  is $<\epsilon/2$.
     Now let $c_1: \R\ra H_2$
    and $c_2:  \R\ra H_2$  be  parametrizations of $q'_xp'$ and
    $q'_yp'$ respectively such that  $c_1(\i)=c_2(\i)=p'$  and
    $c_1(t)$ and $c_2(t)$ are on the
    same horosphere centered at $p'$.
       There are $t_1, t_2\in \R$ with
         $ x'=c_1(t_1)$  and $y'=c_2(t_2)$.
      We may assume  $t_2\ge t_1$.
      Pick  $w\in  P_{p'q'_x}(\p f(E_{x}))$  such that  $d(y',w)<\epsilon/2$.
         By the definition of $x'$ we have  $w=c_1(t_0)$ for some $t_0\le t_1$.
  Since $c_1(t_2)$ and $c_2(t_2)=y'$   are on the same horosphere
    centered at $p'$
  and  $c_1(t_2)$  is the point on the horosphere  that is closest
  to $w$, we have $d(w,  c_1(t_2))\le d(w, y')< \epsilon/2$.
    Consequently $d(x', w)<\epsilon/2$ and $d(x', y')<\epsilon$.

We  next prove the claim.
  Below we shall define $\rho_i=\rho_i(\epsilon)>0$ ($i=1,2,3,4$).
    Set  $\delta_1=\min\{\rho_1,  \rho_2,  \rho_3,  \rho_4\}$.
  Assume  $d_{x'}(q'_x, q'_y)<\delta_1$  and
    $HD_{d_{x'}}(\p f(E_x), \p f(E_y))<\delta_1$.
Let  $z\in P_{p'q'_y}(\p f(E_{y}))$.  Let  $\xi'\in \p f(E_y)$
   such
that its projection on $p'q'_y$ is $z$.   Pick
  some $\eta'\in \p
f(E_x)$ with $d_{x'}(\eta',
        \xi')<\delta_1$.    Let $w$ and $w'$  respectively  be the projections of $\eta'$
          and $\xi'$
        on $p'q'_x$.
   We shall prove
   $d(w, w')<\epsilon/4$  and
    $d(w', z)< \epsilon/4$.

We first show that there exists $\rho_1=\rho_1(\epsilon)>0$ with the
following property:  for any $u, v\in H_1$,
   any  $\xi'\in  \p f(E_v)$,   and any  $\eta'\in \p f(E_u)$,
     if  $d_{u'}(\xi', \eta')<\rho_1$,  then
    $d(w, w')<\epsilon/4$, where $w$  and $w'$  are respectively  the projections of $\eta'$
  and $\xi'$ on $p'q'_u$.
Let $\epsilon'_1=\epsilon'_1(\epsilon/4)$  be given by  Lemma
\ref{l2.2}.
  Since $P_{p'q'_u}(\p f(E_{u}))$   has length at most $B_5$,  Lemma
  \ref{l2.5}
  and the property of visual metrics
  imply   that there is some $\rho_1=\rho_1(\epsilon'_1)>0$ such that
  for any $s\in  P_{p'q'_u}(\p f(E_{u}))$ and $\xi_1, \xi_2\in \p H_2$, if
    $d_{u'}(\xi_1, \xi_2)<\rho_1$, then
      $\angle_s(\xi_1, \xi_2)<\epsilon'_1/10$.
    Let   $\xi'\in  \p f(E_v)$,     $\eta'\in \p f(E_u)$, and
       assume    $d_{u'}(\xi', \eta')<\rho_1$.
       Let $w$ and $w'$  respectively  be the projections of $\eta'$
          and $\xi'$
        on $p'q'_u$.  Then $\angle_w(\xi', \eta')<\epsilon'_1/10$.
 Hence $\angle_w(w', \xi')> \pi/2-\epsilon'_1/10$. Since
 $\angle_{w'}(w, \xi')=\pi/2$, Lemma \ref{l2.2}  implies
 $d(w,w')<\epsilon/4$.

It  remains to prove $d(w', z)<\epsilon/4$.


 Let $\epsilon'_1$ be  as above
and $\epsilon''_1=\epsilon''_1(\epsilon/4, \epsilon'_1)$
  be given by Lemma
\ref{l2.3}.
  Let $c_1$, $c_2$   and  $t_1$,  $t_2$  be as above,
    and $t_0\in \R$  such that
      $c_1([t_0, t_1])=P_{p'q'_x}(\p f(E_{x}))$.
         Then $t_1-t_0\le B_5$.
Lemma \ref{l2.6} implies that there is some
$\rho_2=\rho_2(\epsilon''_1)>0$ such that if
  $d_{x'}(q'_x,q'_y)<\rho_2$,    then
      $d(c_1(t), c_2(t))<\epsilon''_1$  for all $t\in [t_0-\epsilon,
        t_1+\epsilon]$.

         Notice
        $w'=c_1(t')$  for some  $t'\in  [t_0-\epsilon/4,
        t_1+\epsilon/4]$.
  Also $z=c_2(t'')$ for some $t''\in \R$.
We  first assume $t''\ge t_0-\epsilon$.
  The choice of $\rho_2$
      and the convexity of distance function imply
 $d(z, c_1(t''))<\epsilon''_1$ and
   $d(w', c_2(t'))<\epsilon''_1$.
     Suppose
       $d(z, w')\ge \epsilon/4$.
         Then   Lemma \ref{l2.3}  implies
      $\angle_{w'}(c_1(t''), z)<\epsilon'_1/10$. Hence
         $\angle_{w'}(z, \xi')\ge \pi/2-\epsilon'_1/10$.
           Similarly,  $\angle_{z}(w', \xi')\ge\pi/2-\epsilon'_1/10$.
             It follows that
               $\angle_{w'}(z, \xi')+\angle_{z}(w', \xi')\ge
               \pi-\epsilon'_1/5$, contradicting
                 Lemma \ref{l2.2}.

Now   assume $t''\le t_0-\epsilon$.
  Let $\rho_3=B_3/4$, where $B_3$ is the constant in Lemma
  \ref{l3.3}.
  Then
  $$d_{x'}(q'_y, \xi')\ge  d_{x'}(q'_x, \eta')-
  d_{x'}(q'_x, q'_y)-d_{x'}(\eta', \xi')\ge
    {B_3}/2.$$
Hence $d(x', \xi'q'_y)\le D$ for some constant $D=D(L,A)$.
   Noting
  $d(w', x')\le d(w',  w)+d(w, x')\le
    \epsilon/4+B_5$, we have
          $d(w',\xi'q'_y)\le D+B_5+\epsilon/4$.
  Since $d(w', c_2(t'))\le \epsilon''_1$ and
    $c_2(t')\in  zp'$,  Lemma \ref{l2.4} implies that
      $w'$ is a $D'$-quasicenter
        of $p', \xi', q'_y$, where
        $D'=D+B_5+\epsilon/4+C_1+\epsilon''_1$.
  Meanwhile, $z$ is also a $C_1$-quasicenter of
$p', \xi', q'_y$.
  It follows that
  $d(z, w')\le D''$  for some $D''=D''(\epsilon)$.
Lemma \ref{l2.6} implies that there is some
$\rho_4=\rho_4(\epsilon)>0$ such that if
  $d_{x'}(q'_x,q'_y)<\rho_4$,    then
      $d(c_1(t), c_2(t))<\epsilon''_1$
 for all $t\in [t_0-D''-\epsilon,
      t_1+D''+\epsilon]$.  Then the argument in the preceding paragraph
         shows $d(w', z)<\epsilon/4$.

\end{proof}

\b{Le}\label{l4.2} { Given any $\epsilon_1>0$,  there is some
  $\epsilon_2=\epsilon_2(\epsilon_1)>0$  with the following property:
     for any $x,y\in H_1$,  if
        $d_x(q_x, q_y)<\epsilon_2$ and
         $HD_{d_x}(E_x, E_y)< \epsilon_2$,  then
$d_{x'}(q'_x, q'_y)<\epsilon_1$  and
    $HD_{d_{x'}}(\p f(E_x), \p f(E_y))<\epsilon_1$.

}

\end{Le}

\b{proof} For $\epsilon_1>0$, let
  $\epsilon_2$ be determined by
  $\eta(\epsilon_2/{B_1})B_4=\epsilon_1/2$.  Then
   $\epsilon_2=\epsilon_2(\epsilon_1)$.
     Now suppose
  $d_x(q_x, q_y)<\epsilon_2$.  By Lemma \ref{l3.1}
$\partial f: (\p H_1, d_x)\ra (\p H_2, d_{x'})$  is an
  $\eta$-quasisymmetry.
    For any $\xi\in E_x$, we have
     $$\frac{d_{x'}(q'_y,  q'_x)}{d_{x'}(\xi', q'_x)}\le
     \eta\left(\frac{d_x(q_y, q_x)}{d_x(\xi, q_x)}\right)\le
     \eta\left(\frac{\epsilon_2}{B_1}\right).$$
       It follows that $d_{x'}(q'_y,  q'_x)\le
\eta(\frac{\epsilon_2}{B_1})d_{x'}(\xi', q'_x)\le B_4
\eta(\frac{\epsilon_2}{B_1})=\epsilon_1/2<\epsilon_1$.
  The second inequality is
    proved   similarly.

\end{proof}

The following lemma follows from Lemma \ref{l2.6} (2).

\b{Le}\label{l4.3} { Given any $\epsilon_2>0$, there is some
$\delta=\delta(\epsilon_2)>0$ with the following property:
  for any $x, y\in H_1$,   if $d(x,y)<\delta$,  then
    $d_x(q_x, q_y)<\epsilon_2$.

}

\end{Le}

\b{Le}\label{l4.4} { Given any $\epsilon_2>0$, there is some
$\delta=\delta(\epsilon_2)>0$ with the following property:
  for any $x, y\in H_1$,   if $d(x,y)<\delta$,  then
 $HD_{d_x}(E_x, E_y)< \epsilon_2$.

}

\end{Le}

\b{proof} Fix $\epsilon_2>0$.
  Set $\epsilon'_2=\frac{\epsilon_2}{eC_0^2}$.
We shall show that there is some $\delta=\delta(\epsilon'_2)$,
 $1>\delta>0$ with the following
property:
  for any $x, y\in H_1$,   if $d(x,y)<\delta$,
   then for any $\xi\in   E_y$, there is some $\eta\in E_x$ with
   $d_x(\eta, \xi)<\epsilon'_2$.  By symmetry, for any $\beta_1\in E_x$,
   there is some $\beta_2\in E_y$ with
     $d_y(\beta_1, \beta_2)<\epsilon'_2$.
   Since  $\delta<1$,
 the property of visual metrics implies that
$d_x(\beta_1, \beta_2)<\epsilon_2$.

By Lemma \ref{l2.5}, there is some
$\delta_1=\delta_1(\epsilon'_2)>0$ such that if $x\in H_1$ and
$\xi,\eta\in \p H_1$ with $\angle_x(\xi,\eta)< \delta_1$, then
$d_x(\xi, \eta)<\epsilon'_2$. By Lemma \ref{l2.7},  there is some
$\delta_2=\delta_2(\delta_1)>0$ with the following property:
  for any $y_1, y_2\in H_1$ and $\xi\in \p H_1$,  if
    $d(y_1, y_2)<\delta_2$  and $\angle_{y_1}(y_2, \xi)\le \pi/2$,
       $\angle_{y_2}(y_1, \xi)<\pi/2$, then $\angle_{y_2}(y_1,
       \xi)> \pi/2-\delta_1$.
    Let $\xi\in E_y$  and       $w$ its  projection on
      $pq_x$.


\noindent
      {\bf{Claim:}}
 there is some $\delta=\delta(\delta_2)>0$ such that
       $d(x,w)<\delta_2$   whenever   $d(x,y)<\delta$.

  We first finish the proof  assuming the claim.
    We may assume $w\not=x$, otherwise $\xi\in E_x$ and we may
    choose $\eta=\xi$.   Then $\angle_w(x, \xi)=\pi/2$  and $\angle_x(w,
      \xi)<\pi/2$.
       The choice of $\delta_2$  and the claim
      imply   $\angle_x(w, \xi)> \pi/2-\delta_1$.   Let $v\in T_x H_1$ be a
      vector  perpendicular to
      $pq_x$  such that  the angle  between $v$ and the direction of
      $x\xi$ is  $<\delta_1$.  The geodesic ray starting from $x$ in
      the direction of $v$ determines a point $\eta\in E_x$.
        Now the choice of $\delta_1$ implies $d_x(\eta, \xi)<\epsilon'_2$.

We  next   prove the claim.
   Below we shall define $\delta'=\delta'(\delta_2)>0$,
  $\delta''=\delta''(\delta_2)>0$  and
  $\delta'''=\delta'''(L, A)>0$.  Set $\delta=\min\{\delta',
  \delta'', \delta''', 1/2\}$.
   By Lemma \ref{l2.2}, there is a
constant $\epsilon'_1=\epsilon'_1({\delta_2}/2)>0$ with the
following property:
   for any $y_1, y_2, y_3\in H_1$, if $d(y_1, y_2)\ge \delta_2/2$,
   and $\angle_{y_1}(y_2, y_3), \angle_{y_2}(y_1, y_3)\ge 3\pi/8$,
   then $\angle_{y_1}(y_2, y_3)+ \angle_{y_2}(y_1,
   y_3)\le \pi-\epsilon'_1$.
Let $\epsilon_1''=\epsilon_1''(\delta_2/2,  \epsilon'_1)$ be the
constant in Lemma \ref{l2.3}.
  Set  $\delta'=\min\{\epsilon_1'', \delta_2/10\}$.
    By Lemma \ref{l4.3},  there is some $\delta'''=\delta'''(L,A)>0$
      such that $d_y(q_x, q_y)<B_3/2$  whenever $d(x,y)<\delta'''$.
     Assume $d(x,y)<\delta$. To prove $d(x, w)<\delta_2$
        it suffices to show
      $d(y, w)<\delta_2/2$.  We suppose
$d(y, w)\ge \delta_2/2$  and will get a contradiction.  We use the
argument  in the proof of Lemma \ref{l4.1}.
  There are two cases, depending on whether $w$ lies between $x$ and
  $p$  on $pq_x$. If $w$ lies between $x$ and $p$, then
$d(w, yp)<\delta$, and the argument shows
  $\angle_w(y, \xi), \angle_y(w,\xi)\ge \pi/2-\epsilon'_1/10$,
  contradicting Lemma \ref{l2.2}.
 If $x$ lies between $w$ and $p$, then a similar  argument
   as  in the proof of
 Lemma \ref{l4.1} shows that
     $d(w, y)\le D$  for some  constant $D=D(L, A)$.
      Then  by Lemma \ref{l2.6}
        there is some $\delta''=\delta''(\epsilon''_1)>0$ such that
    $d(y, pq_x), d(w, pq_y)<\epsilon''_1$  whenever
$d(x,y)<\delta''$.
  Then one gets a contradiction as above.


\end{proof}

\subsection{$F^{-1}$ is uniformly continuous}\label{s5}

In this Section we  prove that $F^{-1}$  is also uniformly
continuous. Notice that  metric
    spheres in $H_i$ are separating.
Since $F$ is a homeomorphism,   the following proposition implies
that $F^{-1}$ is uniformly continuous:

\b{Prop}\label{p5.1} { For any $0<\epsilon<1$,   there exists a
constant $\delta=\delta(\epsilon)>0$ with the following property:
  for any $x\in H_1$,  $d(F(x),  F(S(x,\epsilon)))\ge \delta$.
    Here $S(x, \epsilon)$  denotes the metric sphere  in $H_1$ with center $x$  and radius $\epsilon$.}

\end{Prop}

 Proposition \ref{p5.1}
   follows from  Lemmas \ref{l5.1}--\ref{l5.3}.

\b{Le}\label{l5.1} { Given any $\epsilon>0$,
 there is some $\delta'=\delta'(\epsilon)>0$ with the following property:
for any $x, y\in H_1$ and $p\in \p H_1$, if $y\in xp$ and $d(x,y)\ge
\epsilon/2$,  then   $d(x', y')>2 \delta'$.

}

\end{Le}

\b{proof}  Assume
  $y\in xp$ and $d(x,y)\ge
\epsilon/2$.
 We first show that there is some $\delta_1=\delta_1(\epsilon)>0$ such that
  $d_x(\xi_1, \xi_2)\ge \delta_1$ for all
  $\xi_1\in E_y$ and $\xi_2\in E_x$.
Let $\epsilon'_1=\epsilon'_1(\epsilon/2)$ be the constant in Lemma
\ref{l2.2}.
   Fix any $\xi_1\in E_y$, $\xi_2\in E_x$.  The assumption
    $d(x,y)\ge \epsilon/2$ and Lemma \ref{l2.2} imply
      $\angle_x(y, \xi_1)\le \pi/2-\epsilon'_1$.
        Hence $\angle_x(\xi_1,\xi_2)\ge \epsilon'_1$.
          Lemma \ref{l2.5}(1) implies that
             there is some $\delta_1=\delta_1(\epsilon'_1)>0$
          such that $d_x(\xi_1,\xi_2)\ge \delta_1$.

Since $\p f:( \p H_1, d_x)\ra (\p H_2, d_{x'})$  is an
 $\eta$-quasisymmetric map, for any $\xi_1\in E_y$, $\xi_2\in E_x$, we
 have
 $$\frac{B_3}{d_{x'}(\xi'_1, \xi'_2)}\le
 \frac{d_{x'}(q'_x, \xi'_2)}{d_{x'}(\xi'_1, \xi'_2)}\le
 \eta\left(\frac{d_x(q_x, \xi_2)}{d_x(\xi_1, \xi_2)}\right)\le
 \eta\left(\frac{B_2}{\delta_1}\right).$$
  It follows that
    $d_{x'}(\p f(E_x), \p f(E_y))\ge  \delta_2=\delta_2(\epsilon)$,
  where  $\delta_2=
     B_3(\eta(\frac{B_2}{\delta_1}))^{-1}$.

By Lemma \ref{l2.5}(2), there is a constant
$\delta_3=\delta_3({\delta_2}/{eC_0^2})>0$:
  for any $z\in H_2$ and $\eta_1, \eta_2\in \p H_2$, if
  $\angle_z(\eta_1, \eta_2)<\delta_3$, then
    $d_z(\eta_1, \eta_2)<\frac{\delta_2}{eC_0^2}$.
    Let $\delta_4=\delta(\delta_3, \lambda)$ be the constant in Lemma
        \ref{l2.7}.  We shall show that $\delta':=\frac{1}{2}\min\{\delta_4,1\}$
     has  the required property. 

      Suppose  $d(x', y')\le  2 \delta'\le \delta_4$.
        Fix some  $\eta'\in \p f(E_x)$ such that the projection of
        $\eta'$   on  $p'q'_x$  is $x'$.
  Lemma \ref{l2.7} applied to $x'$, $y'$ and $\eta'$ implies
$\angle_{y'}(x', \eta')> \pi/2-\delta_3$.
  Denote by $v_1\in T_{y'}H_2$ the initial direction of
    $y'\eta'$ at $y'$.  Let $v_2$ be the projection of $v_1$
   on the hyperplane of $T_{y'}H_2$ orthogonal to the direction of
     $p'q'_x$,  and
 $v_3$   the unit vector in the direction of
                $v_2$.
  Then  the angle between
             $v_1$ and $v_3$ is  less than  $\delta_3$.
         The geodesic in  the sphere  $S_{y'}H_2$
          from $v_1$ to $v_3$  gives rise to a continuous path
             $\gamma $ in $\p H_2$   from $\eta'\in \p f(E_x)$ to a
             point  $\beta$ in $E_{y'}$.  Since the angle between
             $v_1$ and $v_3$ is
  less than  $\delta_3$,  for any $\xi\in \gamma$,
  we have $\angle_{y'}(\eta', \xi)< \delta_3$.
  Now by the choice of $\delta_3$ we have
   $d_{y'}(\xi, \eta')<{\delta_2}/({eC_0^2})$ for all $\xi\in\gamma$.
    Since  $d(x', y')\le  1$, the property of visual
    metric implies  $d_{x'}(\xi, \eta')< \delta_2$.
Since $ \eta'\in \p f(E_x)$  and
 $d_{x'}(\p f(E_x), \p f(E_y))\ge  \delta_2$,
   the path $\gamma$   must  
      lie    in
the component  $B'$  of $\p H_2\backslash \p f(E_y)$ that contains
 $\p f(E_x)$.  In particular,  $\beta\in B'$.
    Notice $B'=\partial f(B)$, where $B$ is the component
      of $\partial H_1\backslash  E_y$ that contains $E_x$.
         Hence $\beta= \xi'_3$ for some $\xi_3\in B$.
There is some $z\in yq_x$, $z\not=y$ such that $\xi_3\in E_z$.
  Since $F$ is a homeomorphism, $z'\in y'q'_x$, $z'\not=y'$.
   It follows that  the projection of
     $\beta$ on $p'q'_x$  lies in $z'q'_x$,
      contradicting
$\beta\in E_{y'}$.


\end{proof}

\b{Le}\label{l5.2} { There is a constant
$\epsilon_2=\epsilon_2(\epsilon)>0$
  with the following property:
     for any $x, y\in H_1$,  any $p\in \p H_1$,  if
     $d(x,y)=\epsilon$ and $d(x, pq_y)<\epsilon_2$,  then
       $d(x', y')\ge \delta'$,
where  $\delta'$ is the constant in  Lemma \ref{l5.1}. }

\end{Le}

\b{proof} Since    $F$ is uniformly continuous, there
  exists   some
$\epsilon_1=\epsilon_1(\delta')>0$ such that
   $d(y'_1, y'_2)<\delta'$ for all $y_1, y_2\in H_1$
       satisfying   $d(y_1, y_2)<\epsilon_1$,
   where  $\delta'$ is the constant in Lemma \ref{l5.1}.
Set $\epsilon_2=\min\{\epsilon/10, \epsilon_1\}$. Now let $x, y\in
H_1$  and  $p\in \p H_1$  with
     $d(x,y)=\epsilon$ and $d(x, pq_y)<\epsilon_2$.
  Let $z\in pq_y$ with $d(x, z)<\epsilon_2$.
  Then  $d(x', z')<\delta'$.
  On the other hand,
  $d(y, z)\ge \epsilon/2$. Lemma \ref{l5.1} implies
    $d(z', y')\ge 2 \delta'$.
      It follows that $d(x', y')\ge \delta'$.

\end{proof}

\b{Le}\label{l5.3} {
 There is a constant $\delta''=\delta''(\epsilon)$ with the
 following property:
 for any $x, y\in H_1$,  any $p\in \p H_1$,  if
     $d(x,y)=\epsilon$ and $d(x, pq_y)\ge \epsilon_2$,  then
       $d(x', y')\ge \delta''$,  where $\epsilon_2$ is the constant
       in Lemma \ref{l5.2}.

}

\end{Le}

\b{proof} Since $d(x, pq_y)\ge \epsilon_2$,
  Lemma \ref{l2.6}(1)  implies
    $d_x(q_x, q_y)\ge \delta_2$, where
    $\delta_2$
     depends only on $\epsilon_2$.
Then the argument in  the
  second paragraph of the proof of Lemma
\ref{l5.1} shows that
  $d_{x'}(q'_x, q'_y)\ge \delta_3$, where $\delta_3$ depends only on
  $\epsilon$.
  Now Lemma \ref{l2.6}(2) implies that
    $d(x', p'q'_y)\ge \delta''$  for some constant
      $\delta''=\delta''(\delta_3)$.
        Since $y'\in p'q'_y$,  the lemma follows.

\end{proof}

 Now   Proposition \ref{p5.1}   holds with  $\delta=\min\{\delta',
 \delta''\}$.

\subsection{Modifying $F$}\label{s6}

In this Section we modify $F$ to obtain a bilipschitz homeomorphism.
  The argument is  a minor modification  of
    that in  Section 7 of \cite{TV2}.
    We  include  it here mainly for completeness.
  Hidden behind this is Sullivan's theory of
  Lipschitz structures.  It is used in the proof of
Lemma  3.9 in \cite{TV2}, which we record as Lemma \ref{l6.3} below.

 Let   $f: X\ra Y$ and $g:X\ra Y$ be two maps
   between metric spaces.
   The distance between $f$ and $g$
    is $d(f, g):=\sup\{d(f(x), g(x)): x\in X\}$.
   We say $f$ and $g$ lie at finite distance from each other if
    $d(f,g)<\i$.
      For any subset $A\subset X$, we also denote
        $d(f,g; A):=d(f|_A, g|_A)$.


Let $f:H_1\ra H_2$  be a quasiisometry and
  $F:  H_1\ra H_2$ the homeomorphism constructed in Section
  \ref{s3}.

\b{Prop}\label{p6.1} {Suppose  $H_1$ and $H_2$ are not
 $4$-dimensional.  Then
   for
  any $\epsilon>0$,  there is a bilipschitz homeomorphism
    $G:  H_1\ra H_2$  such that $d(F, G)\le \epsilon$.
}

\end{Prop}

\noindent
 {\bf{A decomposition of $H_1$}}

  Recall that  $H_1$ is
     the universal cover of a compact Riemannian manifold
${H_1}/{G_1}$
 with sectional curvature $-\lambda^2\le K\le -1$, where $G_1$
    acts on $H_1$ as  the  group  of deck transformations. Fix a compact
   connected
fundamental domain $D_1\subset H_1$ for the action of $G_1$ on
$H_1$.
  Then  $H_1$ is covered by $\{g(D_1):  g\in G_1\}$.
  Let $\Gamma$  be the incidence graph of this covering:
  the vertex set of $\Gamma$ is   $G_1$;
  two vertices $g_1$ and $g_2$ are connected by an edge if
    $g_1(D_1)\cap g_2(D_1)\not=\emptyset$.
  The action of $G_1$ on $H_1$ induces an action   of $G_1$ on   $\Gamma$
  and this action is transitive   on the vertices.  In particular,
     there is some positive integer $m$ such that each vertex has
     valence $m$.

Recall that a graph  is \e{finitely colorable}
     if there is an
integer $N\ge 1$  and a map  $\phi: V\ra \{1,2,\cdots, N\}$ defined
 on  the vertex set  $V$ of the graph, 
such that
  $\phi(v_1)\not=\phi(v_2)$ whenever $v_1$ and $v_2$ are connected by an
  edge.   One   may    assume   $\phi$ is surjective.
     The following well-known  (and easy to prove) lemma
         implies that the  graph  $\Gamma$  defined
  above is  finitely colorable.

\b{Le}\label{l6.1} { Let $\Gamma$ be a graph.  If there is some
integer $m$ such that each vertex has valence at most $m$, then
    $\Gamma$  is finitely colorable.

}

\end{Le}

Hence there is a surjective map
 $\phi:   G_1 \ra  \{1,2,\cdots, N\}$
 such that
    $g_1$ and $g_2$ have different colors
     (i.e., $\phi(g_1)\not=\phi(g_2)$)  whenever
       they are joined by an edge.
  For $1\le i\le N$,
  set
  $${\K}_i=\{g(D_1):  \phi(g)=i\} \;\;  \text{and}\;\;
     {\K}^*_i={\K}_1\cup\cdots \cup {\K}_i.$$
  Notice that
  each family ${\K}_i$  consists of disjoint translates of $D_1$, and
  the $N$ families ${\K}_1, \cdots,  {\K}_N$  cover $H_1$.
Since the action of $G_1$ on $H_1$ is proper and cocompact,
   there is a
    constant $b_1>0$  such that
 $d(g(D_1), h(D_1))\ge b_1$   whenever
   $g(D_1)\cap h(D_1)=\emptyset$.  

\vspace{3mm}

\noindent
  {\bf{A solid family of maps}}

The notion of solid family is defined and discussed in
   \cite{TV2}.  It is closely related to the approximation of
   embeddings  by locally bilipschitz maps.

 Let $X$ and $Y$ be metric spaces. A family $\cal{F}$
  of embeddings $f: X\ra Y$  is said to be
 \e{solid}  if its closure is a compact family of embeddings.
 If $Y=\R^n$  and  $X\subset \R^n$ is either  open or compact,
 then
  $\cal{F}$  is  solid  if and only if
       the following three conditions hold (see p. 315 of \cite{TV2}): \newline
      (1) For every $x_0\in X$, the set $\{f(x_0):  f\in {\cal F} \}$  is
      bounded;\newline
        (2)  For every $x_0\in X$ and $\epsilon>0$,
          there is a neighborhood $U$  of $x_0$
 such that $d(f(x), f(x_0))<\epsilon$  whenever
   $x\in U$ and $f\in \cal{F}$;\newline
     (3)  For every $x_0\in X$  and every neighborhood $U$ of $x_0$,
      there is $\epsilon'>0$ such that
        $d(f(x), f(x_0))\ge \epsilon'$   whenever $x\in X\backslash U$
         and $f\in \cal{F}$.

We  next construct a  solid family from the map $F$.

Let $f:H_1\ra H_2$  be a quasiisometry and
  $F:  H_1\ra H_2$ the homeomorphism constructed in Section
  \ref{s3}.  Let $i=1$ or $2$.  Then
$H_i$ is  the universal cover of a compact Riemannian manifold
${H_i}/{G_i}$  whose sectional curvature satisfies $-\lambda^2\le
K\le -1$.
     Fix a  compact   connected
     fundamental
  domain $D_i$ for $G_i$ and an interior point $x_i$  of $ D_i$.
    Let $d_i$ be the diameter of $D_i$.
 There  is a
    constant $b_i>0$  such that
 $d(g(D_i), h(D_i))\ge b_i$   whenever
   $g(D_i)\cap h(D_i)=\emptyset$.

Let  $A=N_{{b_1}/{2}}(D_1)$  and $A'=N_{{3b_1}/{8}}(D_1)$,
  where  for $t>0$ and a subset $Z\subset X$ of a
    metric space $X$,  $N_t(Z):=\{x\in X:  d(x, Z)<t\}$ denotes the
    open $t$-neighborhood of $Z$.
    Set   $B=\cup\{g(D_1):  g(D_1)\cap D_1\not=\emptyset\}$.
    Then $B$  is compact  connected and $A\subset B$.
  Since $F$ is uniformly continuous, there exists some
   constant $\alpha>0$ such that the diameter of
     $F(g_1(B))$  is at most $\alpha$ for all
  $g_1\in G_1$.

The exponential map   
        $h_i:=\exp_{x_i}: T_{x_i}H_i\ra H_i$ is a diffeomorphism.
           Denote  $U=h_1^{-1}(A)$ and
             $U'=h_1^{-1}(A')$.
 Then   there is some $M_i\ge 1$ such that
 $h_1:  U\ra A$  is a $M_1$-bilipschitz map  and
 the restriction of  $h_2$  on  the  closed ball
  $\overline B(o, \alpha+d_2+1)\subset T_{x_2} H_2$ is
    $M_2$-bilipschitz.


For every  $g_1\in G_1$, we fix some $g'_1\in G_2$ such that
  $g'_1(F(g_1(x_1)))\in D_2$.
    For each $Q=g_1(D_1)$, let
      $F_Q:  U\ra T_{x_2}H_2$ be defined by
$$F_Q:=h^{-1}_2\circ g'_1\circ F\circ g_1\circ   h_1.$$
Notice   that  $g'_1(F(g_1(B)))\subset \overline B(x_2, \alpha+d_2)$
and  $F_Q(U)\subset\overline B(o, \alpha+d_2)\subset T_{x_2} H_2$.
  Set $${\cal{F}}=\{F_Q:   Q=g_1(D_1),\;   g_1\in G_1\}.$$

\b{Le}\label{l6.2} { The family $\cal{F}$  is  solid.

}

\end{Le}

\b{proof}
  It suffices to verify the three conditions above.
    (1) holds because $F_Q(U)\subset \overline B(o, \alpha+d_2)$ for all
      $F_Q\in {\cal{F}}$.
        (2) is true since $F$ is uniformly continuous,  $G_1$, $G_2$
        act as isometries and $h_1|_U$, $h_2|_{\overline B(o,
        \alpha+d_2)}$ are fixed bilipschitz maps. Similarly
          (3) follows from the fact that $F^{-1}$  is uniformly
          continuous.

\end{proof}

\noindent
 {\bf{Proof of Proposition \ref{p6.1}}}

 The following result is key in the proof.
Recall that for $n\not=4$,
 every topological $n$-manifold has a unique
     lipschitz structure
  and
homeomorphisms between
 $n$-manifolds can be approximated by
    locally bilipschitz  homeomorphisms, see \cite{S} or
        Section 4 of \cite{TV2}.

\b{Le}\label{l6.3} \e{(Lemma  3.9 in \cite{TV2})}
  {
Let $U, U', V, W$  be open sets in $\R^n$  such that
 $$W\subset V\subset U, \; \overline{U'}\subset U, \; \overline{W}\cap
 U\subset V,$$
 and $\overline{U'}$   is compact.    Let $\cal{F}$  be a solid
 family of embeddings $g: U\ra \R^n$.  For $n=4$,  we also
   assume that the members of $\cal{F}$  can be  approximated by
     locally bilipschitz
     embeddings.  
       Let
     $\epsilon>0$, $L\ge 1$.  There
       exist   $\delta>0$  and $L'\ge 1$  with
     the following properties:\newline
      Let $h:V\ra \R^n$  be a
        locally  bilipschitz    embedding, and
        $g\in {\cal{F}}$  be such that
         $d(g,h; V)\le \delta$.  Then there is a
            locally bilipschitz  embedding $h':
         U\ra \R^n$  such that\newline
          (1) $d(h', g; U')\le \epsilon$;
          \newline
           (2) $h'=h$  in $W\cap U'$;\newline
            (3) $h'|_{U'}$  is $L'$-bilipschitz  if $h$ is locally
            $L$-bilipschitz.\newline
              Here $\delta$  depends only  on  $\tau=(U, U', V, W,
              {\cal{F}}, \epsilon)$  and $L'$ depends only on $\tau$
              and $L$.

}

\end{Le}

Let $N$,  ${\K}_i$, ${\K}^*_i$ ($1\le i\le N$)
  and $U$, $U'$     be as above.
   For $Q\in  {\K}^*_N$  and  $1\le i\le N$,  let
   $Q_i=N_{2^{-i-1}b_1}(Q)$.
  Set
$V_i=\cup\{Q_i:  Q\in {\K}^*_i\}$, $W_i=\cup\{Q_{i+1}:  Q\in
{\K}^*_i\}$.
  We also let $V_0=\emptyset$,  $W_0=\emptyset$.

For $Q=g_1(D_1)\in {\K}_i$ and $0<t\le {b_1}/2$, let
 $$V_Q(t)=h_1^{-1}(N_t(D_1)\cap g_1^{-1}(V_{i-1})),
 $$
$$W_Q(t)=h_1^{-1}(N_t(D_1)\cap g_1^{-1}(W_{i-1})).$$
  Note  $W_Q(t)\subset V_Q(t)\subset U$.

Since  $G_1$ acts transitively on the family
  ${\K}^*_N$,  the choice of $b_1$  implies the following:
   for any $0<t\le {b_1}/2$,  there is a finite family
   $S(t)$  such that $V_Q(t)$  and $W_Q(t)$
      belong to  $S(t)$  for every $Q\in {\K}^*_N$.   We apply Lemma
      \ref{l6.3}  with $U$,
        $U'$  above,
           $V=V_Q({b_1}/2)$, $W=W_Q({b_1}/2)$, and
           ${\cal{F}}=\{F_Q\}$.   Since the family $S({b_1}/2)$ is finite,
           we obtain:

\b{Le}\label{l6.4} {
  Let $\epsilon>0$ and $L\ge 1$.  Then there are positive numbers
    $\delta=\delta(\epsilon )\le \epsilon$  and $L'=L'(\epsilon,
    L)\ge L$   with the following property:\newline
       Let $Q\in  {\K}^*_N$,  $h: V_Q({b_1}/2)\ra T_{x_2}H_2$
 be a locally $L$-bilipschitz embedding, and
   $g\in \cal{F}$   be such that
     $d(g,h; V_Q({b_1}/2))\le \delta$.  Then
       there is a
         locally  bilipschitz   embedding $h': U\ra  T_{x_2}H_2$  such that\newline
         (1)   $d(h', g;   U')\le
         \epsilon$;\newline
           (2)  $h'=h$  in $W_Q(3b_1/8)$;\newline
            (3)  $h'|_{U'}$  is
            $L'$-bilipschitz.

}

\end{Le}


Since   $F^{-1}$ is uniformly  continuous,
 there is a number $q>0$  such that
  $d(F(x), F(y))\ge q$  whenever   
  $x,y\in   H_1$  with
     $d(x,y)\ge {b_1}/8$.
    Define $\delta_N\ge \delta_{N-1}\ge \cdots\ge \delta_0>0$   by
      $\delta_N=\min\{q/3, \epsilon, 1\}$  and
        $\delta_{j-1}=\delta({\delta_j}/{M_2})/{M_2}$, where $\delta(\cdot)$ is the
        function in Lemma \ref{l6.4}.
          We also  define
          numbers $L_0\le \cdots \le L_N$ by
            $L_0=1$  and $L_j=M_1M_2L'(\delta_j/{M_2}, M_1M_2L_{j-1})$, where $L'(\cdot, \cdot)$ is
            the function in Lemma \ref{l6.4}.    Observe that the
            sequences
              $(\delta_0, \cdots, \delta_N)$ and $(L_0, \cdots,
              L_N)$  depend only on $\epsilon$.
                We show by induction that the following lemma is
                 true for every integer
                  $j\in[0, N]$:

\b{Le}\label{l6.5} {{\bf{Statement}}  ${\bf{S(j):}}$
 There is an embedding   $F_j: V_j\ra H_2$  with the following
 property:\newline
  (1) $d(F_j, F; V_{j})\le \delta_j$;  \newline
    (2)  $F_j(Q_{j})\subset F(N_{{3b_1}/8}(Q))$    for every
  $Q\in {\K}^*_j$;\newline
    (3)   $F_j$   is locally   $L_j$-bilipschitz.

}

\end{Le}

\b{proof} Since $V_0=\emptyset$,  ${\bf{S(0)}}$ is true.  Suppose
that  ${\bf{S(j-1)}}$ is true.   Thus we have an embedding
$F_{j-1}: V_{j-1}\ra H_2$.  We define
  $F_j(x)=F_{j-1}(x)$  for $x\in W_{j-1}$.
  Recall  that  for each $g_1\in G_1$
we fixed  some $g'_1\in G_2$ such that
  $g'_1(F(g_1(x_1)))\in D_2$.
     Let $Q=g_1(D_1)\in {\K}_j$.
    Then   $F_Q=h^{-1}_2\circ g'_1\circ F\circ g_1\circ   h_1.$
    Set
    $$h_Q=h^{-1}_2\circ g'_1\circ F_{j-1}\circ g_1\circ   h_1|_{V_Q({b_1}/2)}.$$
 ${\bf{S(j-1)}}$    implies  that
   $\text{image}(h_Q)\subset \overline B(o, \alpha+d_2+1)\subset
   T_{x_2}H_2$,
  $h_Q$ is locally $M_1M_2L_{j-1}$-bilipschitz,
     and
        $d(h_Q, F_Q; V_Q({b_1}/2))\le M_2\delta_{j-1}$.
    Hence we can apply Lemma \ref{l6.4}  with
         $g=F_Q, h=h_Q$,
      $\epsilon=\delta_j/{M_2}$  and $L=M_1M_2L_{j-1}$.
    We obtain a  locally
      bilipschitz    embedding
    $h'_Q:  U\ra  T_{x_2}H_2$  such that\newline
      (a) $d(h'_Q, F_Q; U')\le {\delta_j}/{M_2}$,\newline
        (b)  $h'_Q=h_Q$  on $W_Q({3b_1}/8)$,\newline
           (c) $h'_Q|_{U'}$ is  $L'(\delta_j/{M_2},
          M_1M_2L_{j-1})$-bilipschitz.\newline
           Setting
           $$F_j=(g'_1)^{-1}\circ h_2\circ h'_Q\circ h^{-1}_1\circ
           g^{-1}_1$$
             in $Q_j$  we obtain a well-defined  map
             $F_j:  V_j\ra H_2$.  We show that $F_j$  satisfies the
             conditions  (1), (2),   (3) and that $F_j$ is injective.

Let $Q\in {\K}^*_j$.  If $Q\in  {\K}^*_{j-1}$, (1) follows from
     ${\bf{S(j-1)}}$.
  If $Q\in {\K}_j$,
  (a)  implies  $h'_Q(U')\subset \overline B(o,
  \alpha+d_2+1)\subset T_{x_2}H_2$  and hence
    $$d(F_j, F; Q_{j})\le M_2 \,d(h'_Q, F_Q;
      U')\le \delta_j.$$

    To prove (2), let again
$Q\in {\K}^*_j$.  If $Q\in  {\K}^*_{j-1}$, (2) follows from
  ${\bf{S(j-1)}}$.
Suppose  $Q\in {\K}_j$.
  Since $\delta_j<q$,
    the choice
    of $q$   and  ${\bf{S(j)}}$ (1)
   imply  (2).


  If  $Q\in {\K}^*_{j-1}$,  $F_j$  is locally $L_j$-bilipschitz in
  $Q_j$  by   ${\bf{S(j-1)}}$ (3).  If  $Q\in {\K}_j$,  then
    (c) implies $F_j|_{Q_j}$ is $L_j$-bilipschitz.
       Hence   $F_j$  is   a   locally  $L_j$-bilipschitz  immersion.  We
       finally  show that $F_j$ is injective.  We know that
         $F_j|_{Q_j}$ is injective for every $Q\in
         {\K}^*_j$.  Moreover,  if
          $Q, R\in{\K}^*_j$  and $Q\cap R=\emptyset$,  then  (2)
          implies that
            $F_j(Q_j)\cap F_j(R_j)=\emptyset.$
              Hence it suffices  to show that $F_j(x)\not=F_j(y)$
              when $j\ge 2$, $x\not=y$, $x\in Q_j$
              and $y\in  R_j$  where
               $Q\in {\K}_j$,   $R\in {\K}^*_{j-1}$  and $Q\cap
               R\not=\emptyset$.
                 The equality
                   $F_j=(g'_1)^{-1}\circ h_2\circ h'_Q\circ h^{-1}_1\circ
           g^{-1}_1$
                     is valid in $N_{{3b_1}/8}(Q)\cap W_{j-1}$.  Hence
                     we  may assume that
                       $y\notin N_{{3b_1}/8}(Q)$.  By the choice of
                       $q$, we have 
                           $d(F(x), F(y))\ge  q$.
                             By (1)  we obtain
                             $$d(F_j(x), F_j(y))\ge   d(F(x),
                             F(y))-d(F_j(x), F(x))-d(F_j(y), F(y))
                             \ge q-2\delta_j\ge q/3>0.$$

\end{proof}

    Now notice $V_N=H_1$.
      Hence $F_N: H_1\ra H_2$ is an embedding with
       $d(F_N, F)\le \delta_N\le \epsilon$.
It follows that $F_N$ is a homeomorphism.
 By (3),  $F_N$ is
   locally $L_N$-bilipschitz.
  Since $H_1$ and $H_2$
  are geodesic metric spaces,
  $F_N$ is $L_N$-bilipschitz.

This completes the proof of Proposition \ref{p6.1}.


\noindent
 {\bf{Proof of Theorem  \ref{main}}}
Let  $f: H_1\ra H_2$ be a quasiisometry,
 $F$ the map constructed in Section \ref{s3}
   and $G$ the map in Proposition \ref{p6.1}.
      Since $G$ is bilipschitz   and
        $d(F, G)<\i$,  $F$ is also a quasiisometry.
          From the construction of $F$ one sees easily that
            $\partial F=\partial f$.
            We shall prove  $d(f, F)<\i$.

Let $L\ge 1$ and $A\ge 0$ be such that both $f$ and $F$ are
$(L,A)$-quasiisometries.
  Let $x\in H_1$.  Pick $p, q, r\in \partial H_1$  such that
   $x\in pq$ and that $xr$ is perpendicular to $pq$.
     By Lemma \ref{l2.4}
          $d(x, pr),  d(x, qr)\le C_1$.
       So  $d(f(x), f(pr)), d(f(x), f(qr))\le L\, C_1+A$.
 Since  $f(pr)$, $f(qr)$ are $(L,A)$-quasigeodesics, the stability
 of quasigeodesics yields
   $HD(f(pr), p'r')\le C$,
$HD(f(qr), q'r')\le C$
  and  $HD(f(pq), p'q')\le C$,
 where $C$ depends only on $L$  and  $A$.
 It follows that
  $d(f(x), p'r')\le C+L\, C_1+A$,
$d(f(x), q'r')\le C+L\, C_1+A$
  and $d(f(x), p'q')\le C$.
    In other words, $f(x)$ is a
      $(C+L\, C_1+A)$-quasicenter of $p', q', r'$.
       Similarly $F(x)$ is also a
$(C+L\, C_1+A)$-quasicenter of $p', q', r'$.
  It follows that $d(F(x), f(x))\le D$ for some $D$ depending only on
  $L$  and  $A$.   This is true for every $x\in H_1$. Hence
   $d(F, f)\le D$.

   \qed

\noindent
 {\bf{Proof of Corollary   \ref{ch}}}
  Let  $h: \partial B^n_\C\ra \partial B^n_\C$  ($n\not=2$)
      be a quasisymmetric
  map,  where  $\partial B^n_\C$  is equipped with
    the  Carnot metric.
    Then  there is a quasiisometry $f: B^n_\C\ra B^n_\C$
    with $\partial f=h$ (see \cite{BS}). By Theorem \ref{main},
      there is a bilipschitz homeomorphism
       $G:B^n_\C\ra B^n_\C$  with
         $d(G, f)<\i$.   $G$ is clearly a quasiconformal map in the
          complex hyperbolic metric.   The fact
$d(G, f)<\i$  implies that $G$ and $f$ have the same boundary map,
which is $h$.

\qed

\subsection{Open   questions}\label{s7}

 In this Section
we present
   several
  questions related to the result in this paper.
 The first natural question is the following.  

\b{question}\label{p1} {Let $H_1$ and $H_2$ be two Hadamard
  $n$-manifolds
  (whose   sectional curvatures  are  bounded from below)
    with $n\not=4$,  and
   $f: H_1\ra H_2$ a quasiisometry. Is $f$ always a finite distance
   from a bilipschitz homeomorphism?

}

\end{question}

Notice that a Hadamard manifold has bounded geometry
  if the sectional curvature is bounded from below.

 Recall that  a subset $A\subset X$ of  a metric space $X$ is a
   \e{separated net}
 if there are constants $a, b>0$ such that
   $d(x,y)\ge a$ for distinct $x, y\in A$ and
      $d(x, A)\le b$ for all $x\in X$.
   Observe  that the restriction of
  a quasiisometry $f:  X\ra Y$ to a suitable
    separated net   is   a bilipschitz
  embedding of the net into  $Y$.
     Hence  an affirmative answer to the following problem implies
     the positive answer to Question \ref{p1}.

\b{question}\label{p2} {Let $H_1$ and $H_2$ be two Hadamard
  $n$-manifolds
  (whose   sectional curvatures  are  bounded from below)
    with $n\not=4$,
  $A\subset H_1$  a  separated net,
   and
   $f: A\ra H_2$   a  bilipschitz embedding. 
  Does  $f$ always  extend to a
    bilipschitz homeomorphism  from $H_1$ to $H_2$?

}

\end{question}

 In the case $H_1=H_2=\R^n$,  Question \ref{p2}
     has been    asked  by
Alestalo-Trotsenko-Vaisala \cite{ATV}.

Hadamard manifolds have no topology:  they are contractible.
  But Question \ref{p1} can also be asked for more general
  manifolds. For example one can consider quasiisometries between
    noncompact  hyperbolic  surfaces.

\b{question}\label{p4} {
 Let $X$ and $Y$ be two  open complete Riemannian
   $n$-manifolds ($n\not=4$).
   Suppose $X$ and $Y$ are Gromov hyperbolic and have bounded
   geometry.
     Let $f: X\ra Y$ be a quasiisometry.
       Suppose there is a homeomorphism $F:\overline X\ra \overline
       Y$  such that $F|_{\partial X}=\partial f$.  Is $f$ always
       at  a  finite distance from a bilipschitz homeomorphism?
    Here $\overline X$  and $\overline Y$ are Gromov
    compactifications of  $X$ and $Y$ respectively.
}

\end{question}

  For an arbitrary quasiisometry
 $f: X\ra Y$,  the boundary map
   $\p f$  in general does not have a homeomorphic extension, let
   alone a   bilipschitz extension.

Every quasisymmetric map $f: X\ra Y$ between two metric spaces
extends to a quasisymmetric map between their completions.
  Hence quasisymmetric maps between Euclidean domains extend to
  quasisymmetric maps between their closures.
    A basic question is to what extent the converse is true, that
    is,   under what conditions, a quasisymmetric map between the
    boundaries of two domains extends to a quasisymmetric map
      between   the
    domains?  The theorem of Tukia-Vaisala shows that
     it is the case when the domains are balls in $\R^n$.
  How about more general domains?
    A necessary condition is that the quasisymmetric map between the
    boundaries must be power quasisymmetric:  quasisymmetric maps
    between connected metric spaces are power quasisymmetric,   and
    domains are connected.

Again,    in general
 a quasisymmetric  map
    between the
boundaries may not  have a homeomorphic extension to the closures.

\b{question}\label{p3} {
  Let $\Omega_1, \Omega_2\subset \R^n$ be two  domains,
    and $f:\partial \Omega_1\ra \partial \Omega_2$
     a  power
     quasisymmetric map.  Suppose there is a homeomorphism
       $F: \overline\Omega_1\ra \overline\Omega_2$ such that
       $F|_{\partial \Omega_1}=f$.  Does $f$ extend to a
       quasisymmetric map from $\Omega_1$ to $\Omega_2$?

}

\end{question}

One may  have to restrict
  attention to the so-called uniform domains.
  Uniform domains are considered nice domains in many
    analysis problems.  And they are Gromov hyperbolic in the
     quasihyperbolic metric.   See \cite{BHK}
       for more details.


 \addcontentsline{toc}{subsection}{References}


\begin{thebibliography}{99}



\bibitem[A]{A}
L.  Ahlfors,
  \e{Extension of quasiconformal mappings from two to three
dimensions,}
   Proc. Nat. Acad. Sci. U.S.A. {\bf{51}} 1964 768--771.



\bibitem[ATV]{ATV}
P.  Alestalo, D.A.  Trotsenko, J. Vaisala,
 \e{The linear extension
property of bi-Lipschitz mappings,} (Russian) Sibirsk. Mat. Zh.
{\bf{44}} (2003), no. 6, 1226--1238; translation in Siberian Math.
J. {\bf{44}} (2003), no. 6, 959--968.


\bibitem[B]{B}
M.  Bourdon,
  \e{Structure conforme au bord et flot geodesique d'un
${\rm CAT}(-1)$-espace,}
    Enseign. Math. (2) {\bf{41}} (1995), no. 1-2,
63--102.




\bibitem[BA]{BA}
A.  Beurling,  L.  Ahlfors,
 \e{The boundary correspondence under quasiconformal
mappings,}   Acta Math. {\bf{96}} (1956), 125--142.






\bibitem[BH]{BH}
 M. Bridson, A.  Haefliger,
 \e{Metric spaces of
 non-positive curvature,}
  Grundlehren der Mathematischen Wissenschaften, {\bf{319}}.
     Springer-Verlag, Berlin, 1999.








\bibitem[BHK]{BHK} M. Bonk, J.  Heinonen, P.  Koskela,
 \e{Uniformizing Gromov hyperbolic spaces,}  Asterisque No. {\bf{270}} (2001).









\bibitem[BS]{BS}
M.  Bonk, O.  Schramm,
   \e{Embeddings of Gromov hyperbolic spaces, }
       Geom.
Funct. Anal.  {\bf{10 }}(2000), no. 2, 266--306.




\bibitem[BW]{BW}
J. Block, S.  Weinberger,
  \e{Large scale homology theories and
geometry,}
    Geometric topology (Athens, GA, 1993), 522--569, AMS/IP
Stud. Adv. Math., 2.1, Amer. Math. Soc., Providence, RI, 1997.







\bibitem[C]{C}
L.  Carleson, \e{The extension problem for quasiconformal mappings,}
Contributions to analysis (a collection of papers dedicated to
Lipman Bers), pp. 39--47. Academic Press, New York, 1974.


\bibitem[CDP]{CDP}
M.    Coornaert, T.   Delzant, A.   Papadopoulos,
 \e{Geometrie et
theorie des groupes,
  Les groupes hyperboliques de Gromov}, Lecture Notes in Mathematics, {\bf{1441}}. Springer-Verlag, Berlin, 1990.


\bibitem[CE]{CE}
J.  Cheeger, D.  Ebin, \e{Comparison theorems in Riemannian
geometry,}
   North-Holland Mathematical Library, Vol. {\bf{9}}, 1975.




\bibitem[CF]{CF}
T.  Chapman, S.   Ferry,
 \e{Approximating homotopy equivalences by
homeomorphisms,}
     Amer. J. Math. {\bf{101}} (1979), no. 3, 583--607.






\bibitem[DFW]{DFW}
 A. Dranishnikov, S.  Ferry, S.  Weinberger,
   \e{Large Riemannian manifolds which are flexible,}
       Ann. of Math. (2) {\bf{157}} (2003), no. 3, 919--938.






\bibitem[DS]{DS}
S.  Donaldson, D.   Sullivan,
 \e{Quasiconformal $4$-manifolds,}
     Acta
Math. {\bf{163}} (1989), no. 3-4, 181--252.




\bibitem[G]{G}
W.   Goldman, \e{Complex hyperbolic geometry,}
 Oxford Mathematical Monographs, 1999.




\bibitem[S]{S}
D.  Sullivan,
  \e{Hyperbolic geometry and homeomorphisms,}
       Geometric
topology (Proc. Georgia Topology Conf., Athens, Ga., 1977), pp.
543--555, Academic Press, New York-London, 1979.



\bibitem[TV]{TV}
P.  Tukia, J.  V\"ais\"al\"a,
   \e{Quasiconformal extension from dimension $n$
to $n+1$,}   Ann. of Math. (2) {\bf{115}} (1982), no. 2, 331--348.




\bibitem[TV2]{TV2}
P.  Tukia, J.  V\"ais\"al\"a, \e{Lipschitz and quasiconformal
approximation and extension,}
    Ann. Acad. Sci. Fenn. Ser. A I Math. {\bf{6}}
(1981), no. 2, 303--342 (1982).






\bibitem[V]{V} J. V\"ais\"al\"a,
\e{The free quasiworld,
 Quasiconformal geometry and dynamics} (Lublin, 1996), 55--118, Banach Center Publ., 48, Polish Acad. Sci., Warsaw, 1999.





\bibitem[W]{W}
K. Whyte, \e{Amenability, bi-Lipschitz equivalence, and the von
Neumann conjecture, }
    Duke Math. J. {\bf{99}} (1999), no. 1, 93--112.



\end{thebibliography}
\end{document}